\documentclass{article}

\usepackage{microtype}
\usepackage{graphicx}
\usepackage{subfigure}
\usepackage{booktabs} 

\usepackage{hyperref}

\usepackage[accepted]{icml2020}

\usepackage{multirow}
\usepackage{amsmath,amsthm,latexsym} 
\usepackage[fixamsmath]{mathtools}
\mathtoolsset{showmanualtags,mathic,centercolon} 
\usepackage{amssymb,amsfonts} 

\newcommand{\R}{\mathbb{R}} 
\newcommand{\bigO}{\mathcal{O}} 

\def\be{\begin{equation}}
\def\ee{\end{equation}}

\renewcommand{\t}[1]{\widetilde{#1}} 

\icmltitlerunning{Scalable Derivative-Free Optimization for Nonlinear Least-Squares Problems}

\begin{document}

\twocolumn[
\icmltitle{Scalable Derivative-Free Optimization for Nonlinear Least-Squares Problems}



\icmlsetsymbol{equal}{*}

\begin{icmlauthorlist}
\icmlauthor{Coralia Cartis}{equal,oxmi,ati} 
\icmlauthor{Tyler Ferguson}{equal,oxmi,oxbdi} 
\icmlauthor{Lindon Roberts}{equal,anu} 
\end{icmlauthorlist}

\icmlaffiliation{oxmi}{Mathematical Institute, University of Oxford, Oxford, Oxfordshire, United Kingdom}
\icmlaffiliation{ati}{The Alan Turing Institute for Data Science, London, UK}
\icmlaffiliation{oxbdi}{Big Data Institute, University of Oxford, Oxford, Oxfordshire, United Kingdom}
\icmlaffiliation{anu}{Mathematical Sciences Institute, Australian National University, Canberra, ACT, Australia}

\icmlcorrespondingauthor{Lindon Roberts}{lindon.roberts@anu.edu.au}

\icmlkeywords{derivative-free optimization, nonlinear least-squares, sketching}

\vskip 0.3in
]



\printAffiliationsAndNotice{\icmlEqualContribution} 

\begin{abstract}
Derivative-free---or zeroth-order---optimization (DFO) has gained recent attention for its ability to solve problems in a variety of application areas, including machine learning, particularly involving objectives which are stochastic and/or expensive to compute.
In this work, we develop a novel model-based DFO method for solving nonlinear least-squares problems.
We improve on state-of-the-art DFO by performing dimensionality reduction in the observational space using sketching methods, avoiding the construction of a full local model.
Our approach has a per-iteration computational cost which is linear in problem dimension in a big data regime, and numerical evidence demonstrates that, compared to existing software, it has dramatically improved runtime performance on overdetermined least-squares problems.
\end{abstract}

\section{Introduction}
Derivative-free optimization (DFO), or zeroth order optimization, refers to optimization when no gradient information for the objective (and/or constraints) is available \cite{Larson2019}, such as when the objective is stochastic and/or computationally expensive.
This type of problem arises across a broad range of application areas \cite{Conn2009,Audet2017}, but has attracted particular  recent attention in the learning community for problems such as black-box attacks \cite{Chen2017,Ughi2019}, hyperparameter tuning \cite{Ghanbari2017,Lakhmiri2019} and reinforcement learning \cite{Mania2018,Choromanski2019}.
A current deficiency of DFO methods is their performance on large-scale problems, which is critical to their utility in machine learning; there have been several recent works aimed at improving the scalability of DFO \cite{Bergou2019,Roberts2019,Porcelli2020,Cristofari2020}.

Here, we consider `model-based' DFO methods for large-scale nonlinear least-squares problems \cite{Zhang2010,Wild2017,Cartis2019a}, an important problem class in machine learning \cite{Cai2019}.
Specifically, we improve the scalability of DFO methods for nonlinear least-squares problems in the `big data' regime (i.e.~fitting to many observations).
We do this by using sketching techniques from randomized numerical linear algebra to avoid ever constructing a full local model, yielding an algorithm with: (1) a per-iteration computational cost with linear dependence on dimension in the big data regime, compared to a quadratic dependence for existing methods; and (2) an order-of-magnitude lower runtime compared to state-of-the-art software, while yielding comparable objective reductions in low-accuracy regimes for large-scale problems.

Sketching techniques have attracted substantial attention for large-scale linear algebra problems such as linear least-squares and low-rank approximation \cite{Halko2011,Mahoney2011,Woodruff2014}. 
However, despite success in the linear least-squares setting \cite{Dahiya2018}, to our knowledge they have not been used specifically for nonlinear least-squares problems in a DFO context.
We note that similar techniques have been applied to derivative-based nonlinear least-squares \cite{Ergen2019}, BFGS \cite{Gower2016}, Newton's method \cite{Gower2019,RoostaKhorasani2019,Berahas2020} and SAGA \cite{Gower2020}.
We also note that (randomized) gradient sampling methods for DFO model construction were analyzed in \cite{Berahas2019}.


\section{DFO for Nonlinear Least-Squares} \label{sec_dfols}
In this section we outline a model-based DFO method for nonlinear least-squares problems, a class of methods well-suited to exploiting problem structure and with strong practical performance.
Our problem here is
\vspace{-0.5em}
\be \min_{x\in\R^d} f(x) := \frac{1}{2}\|r(x)\|_2^2 = \frac{1}{2}\sum_{i=1}^{n} r_i(x)^2, \label{eq_nlls} \ee
where $r(x):=[r_1(x), \ldots, r_n(x)]^T:\R^d\to\R^n$, for $n\geq d$.
Motivated by the classical Gauss-Newton method \cite{Nocedal2006}, at iteration $k$ we build a linear model for $r(x)$ which we hope is accurate close to our iterate $x_k$:
\be r(x_k + s) \approx M_k(s) := r(x_k) + J_k s, \label{eq_r_model} \ee
where $J_k\in\R^{n\times d}$.
To find $J_k$, we maintain a collection of points $\{y_0:=x_k, y_1, \ldots, y_d\}\subset\R^d$ and require that $M_k$ interpolate $r$ at these points (i.e.~\eqref{eq_r_model} holds with equality). 
We thus find $J_k$ by requiring $M_k(y_t - x_k) = r(y_t)$ for all $t=0,\ldots,d$ 
which yields the linear system
\be \begin{bmatrix}(y_1-x_k)^T \\ \vdots \\ (y_d-x_k)^T\end{bmatrix} J_k^T = \begin{bmatrix}(r(y_1)-r(x_k))^T \\ \vdots \\ (r(y_d)-r(x_k))^T\end{bmatrix}. \label{eq_interp_system} \ee
Our linear model \eqref{eq_r_model} naturally gives a local convex quadratic model for the objective, $m_k(s)\approx f(x_k+s)$, given by
\be m_k(s) := \frac{1}{2}\|M_k(s)\|_2^2 = c_k + g_k^T s + \frac{1}{2}s^T H_k s, \label{eq_full_model} \ee
where $c_k:=f(x_k)$, $g_k:=J_k^T r(x_k)$ and $H_k:=J_k^T J_k$.
Provided the interpolation points have good geometry (they are close enough to $x_k$ and the linear system in \eqref{eq_interp_system} is well-conditioned), we can guarantee that $m_k$ is a comparably accurate approximation for $f$ as the corresponding derivative-based model \citep[Lemma 3.3]{Cartis2019a}.

This approximation can then be implemented inside a trust-region framework \cite{Conn2000}, with appropriate updating of the interpolation set, to yield a globally convergent algorithm \cite{Cartis2019a}.
The software DFO-LS \cite{Cartis2018} is an implementation of this approach.

\subsection{Scalability of existing methods}
\begin{figure}[tbh]
    \vskip 0.2in
    \begin{center}
        \centerline{\includegraphics[width=\columnwidth]{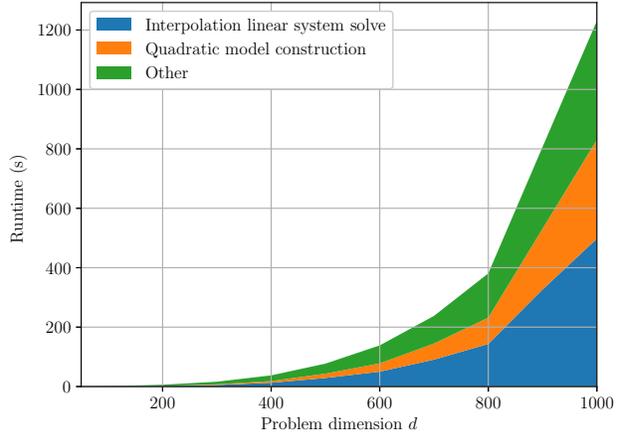}}
        \caption{Runtime of DFO-LS v1.2.1 split by task for the $d$-dimensional generalized Rosenbrock function. The tasks `interpolation linear system solve' and `quadratic model construction' refer to solving \eqref{eq_interp_system} and forming $g_k$ and $H_k$ \eqref{eq_full_model}, respectively.}
        \label{fig_dfols_times}
    \end{center}
    \vskip -0.2in
\end{figure}

In general, model-based DFO methods are best suited to small-scale problems (i.e.~$d$, $n$ small).
An issue impacting their success at scale is the cost of solving the $d\times d$ system \eqref{eq_interp_system} (with $n$ right-hand sides).
To illustrate this, we run DFO-LS on the generalized Rosenbrock function in $\R^d$ (with $n=2d$), defined by
\be r_{2i-1}(x) := 10(x_{i+1} - x_i^2) \quad \text{and} \quad r_{2i}(x) := x_i-1, \ee
where $i=1,\ldots,d$, for a budget of $3(d+1)$ objective evaluations and $d\in[50,1000]$.
In Figure \ref{fig_dfols_times} we show the total runtime of DFO-LS split by the different parts of the algorithm.
As the underlying problem dimension increases, the dominant runtime costs of DFO-LS are: solving the interpolation system \eqref{eq_interp_system}, and evaluating $g_k$ and $H_k$ \eqref{eq_full_model}.

We can quantify the cost of these two steps.
Firstly, solving \eqref{eq_interp_system} requires factorizing a $d\times d$ matrix and backsolving with $n$ right-hand sides, for total cost $\bigO(d^3 + nd^2)$ flops.
In addition, constructing $g_k$ and $H_k$ via standard dense linear algebra requires $\bigO(nd + nd^2) = \bigO(nd^2)$ flops.

Our proposal here is to reduce these costs by the use of sketching techniques from randomized numerical linear algebra for model construction.
This corresponds to dimensionality reduction in $n$, so is particularly useful in the big data regime $n\gg d$, when $r$ is overdetermined.

\section{Sketching for Linear Least-Squares}  \label{sec_sketching}
Sketching is a dimensionality reduction technique that has become a popular tool in numerical linear algebra, including for solving the linear least-squares problem
\be \min_{x\in\R^d} \|Ax-b\|_2^2, \label{eq_lls} \ee
where $A\in\R^{n\times d}$ is a full rank matrix with $n\gg d$.
The idea is to select a random `sketching matrix' $S\in\R^{m\times n}$ with $m\ll n$ and instead solve the smaller $m\times d$ problem
\be \min_{x\in\R^d} \|S(Ax-b)\|_2^2 = \|(SA)x-(Sb)\|_2^2. \label{eq_sketched_lls} \ee
If the distribution for $S$ is chosen well and $m$ is sufficiently large, the minimizer of \eqref{eq_sketched_lls} is close to the minimizer of \eqref{eq_lls} with high probability.
For example, if $S$ is a hashing matrix (see below) and $m\sim d^2/\epsilon$ then solving \eqref{eq_sketched_lls} gives an $\epsilon$-accurate minimizer of \eqref{eq_lls} with probability 0.99 \citep[Theorem 2.16]{Woodruff2014}.
In this case $m$ is \emph{independent of $n$}.

Several choices for $S$ have been proposed, such as Gaussian matrices \cite{Sarlos2006}, subsampling matrices (where each row of $S$ is a randomly-chosen coordinate vector in $\R^n$) \cite{Drineas2006} and hashing matrices (where each column of $S$ has a small number of randomly-chosen nonzero entries, with random value $\pm 1$) \cite{Clarkson2017}.
These methods vary in their requirements on $m$ and the cost of constructing $SA$ and $Sb$.


\section{Sketching in DFO for Nonlinear Least-Squares} \label{sec_dfols_sketching}

We now present our sketching-based DFO algorithm for \eqref{eq_nlls}.
Motivated by \eqref{eq_sketched_lls}, we replace $\|M_k(s)\|_2^2$ in \eqref{eq_full_model} with $\|S_k M_k(s)\|_2^2$ for some sketching matrix $S_k\in\R^{m\times n}$ (resampled for each $k$).
We achieve this by solving the reduced interpolation system (c.f.~\eqref{eq_interp_system})
\be \begin{bmatrix}(y_1-x_k)^T \\ \vdots \\ (y_d-x_k)^T\end{bmatrix} (S_k J_k)^T = \underbrace{\begin{bmatrix}(r(y_1)-r(x_k))^T \\ \vdots \\ (r(y_d)-r(x_k))^T\end{bmatrix}}_{=:R_k} S_k^T, \label{eq_interp_system_sketched} \ee
to find the sketched model Jacobian $(S_k J_k)\in\R^{m\times d}$.
Our new quadratic model for $f$ (c.f.~\eqref{eq_full_model}) is
\be \t{m}_k(s) := \frac{1}{2}\|S_k r(x_k) + S_k J_k s\|_2^2 = \t{c}_k + \t{g}_k^T s + \frac{1}{2}s^T \t{H}_k s, \label{eq_full_model_sketched} \ee
where $\t{c}_k:= \frac{1}{2}\|S_k r(x_k)\|_2^2$, $\t{g}_k:=(S_k J_k)^T S_k r(x_k)$ and $\t{H}_k:=(S_k J_k)^T (S_k J_k)$.
We note that the Jacobian $J_k$ and model $M_k$ are never constructed; we only ever explicitly form $(S_k J_k)$ and $\t{m}_k$.
This model is used inside a trust region method, presented in Algorithm \ref{alg_dfols_sketching}.
Details on updating the interpolation set can be found in \cite{Cartis2018}.

\subsection{Per-iteration computational work of Algorithm \ref{alg_dfols_sketching}} \label{sec_sketch_cost}

The cost of solving the sketched system \eqref{eq_interp_system_sketched} is $\bigO(d^3)$ to factorize the matrix plus $\bigO(md^2)$ to solve with $m$ right-hand sides.
Then we have to form $\t{g}_k$ and $\t{H}_k$, with cost $\bigO(md + md^2)=\bigO(md^2)$.\footnote{We never form $\t{c}_k$ as it does not affect the step $s_k$ or $\rho_k$ \eqref{eq_rho}.}
This reduces the cost of constructing $S_k J_k$ and $\t{m}_k$ so it now depends on $m\ll n$, however we still need to form $R_k S_k^T$ to solve \eqref{eq_interp_system_sketched}. 
Regardless of our choice of $S_k$, we first form the matrix $R_k\in\R^{d\times n}$, with cost $\bigO(nd)$.
Since $R_k$ is typically dense, we have costs:
\begin{itemize}
    \item If $S_k$ is Gaussian, $\bigO(mn)$ to form $S_k$ and $\bigO(mnd)$ to compute $R_k S_k^T$;
    \item If $S_k$ is formed by sampling, then $R_k S_k^T$ is a sample of $m$ columns from $R_k$, so generating $S_k$ costs $\bigO(m)$ and computing $R_k S_k^T$ requires no flops;
    \item If $S_k$ is formed by hashing with $s$ nonzeros per column, then forming $S_k$ costs $\bigO(sn)$ to generate the nonzero indices and values, and building $R_k S_k^T$ costs $\bigO(d\cdot\operatorname{nnz}(S_k))=\bigO(snd)$. 
    Typically we have $s=1$ or 2.
\end{itemize}
\vspace{-0.7em}
That is, the cost of computing $R_k S_k^T$ is $\bigO(nd)$ for sampling and hashing or $\bigO(mnd)$ for Gaussian sketching.
Similarly, forming $S_k r(x_k)$ takes $\bigO(m)$ flops for sampling, $\bigO(n)$ for hashing or $\bigO(mn)$ for Gaussian sketching.


\begin{algorithm}[tb]
    \caption{Model-based DFO for \eqref{eq_nlls} with Sketching}
    \label{alg_dfols_sketching}
\begin{algorithmic}[1]
    \STATE {\bfseries Input:} $r:\R^d\to\R^n$, $x_0\in\R^d$, $0<\Delta_0<\Delta_{\max}$, $0<\eta_1<\eta_2<1$ and $0<\gamma_{\rm dec}<1<\gamma_{\rm inc}$.
    \STATE Evaluate $r$ at $x_0$ and an arbitrary interpolation set.
    \FOR{$k=0,1,\ldots$}
        \STATE Find $S_k J_k$ by solving \eqref{eq_interp_system_sketched} and build $\t{m}_k$ \eqref{eq_full_model_sketched}
        \STATE Approximately solve the trust-region subproblem $\min_{s\in\R^d} \t{m}_k(s) \: \text{s.t.}\: \|s\|_2\leq\Delta_k$ to get a step $s_k$.
        \STATE Evaluate $r(x_k+s_k)$ and calculate
        \be \rho_k := \frac{f(x_k) - f(x_k+s_k)}{\t{m}_k(0)-\t{m}_k(s_k)}. \label{eq_rho} \ee
        \STATE Set $x_{k+1} = x_k + s_k$ if $\rho_k\geq\eta_1$, else $x_{k+1}=x_k$.
        \STATE Set $\Delta_{k+1}=\min(\gamma_{\rm inc}\Delta_k, \Delta_{\max})$ if $\rho_k\geq \eta_2$, $\Delta_{k+1}=\Delta_k$ if $\rho_k\in[\eta_1,\eta_2)$, or $\Delta_{k+1}=\gamma_{\rm dec}\Delta_k$ otherwise.
        \STATE Update the interpolation set to include $x_k+s_k$ and ensure the set has sufficiently good geometry.
    \ENDFOR
\end{algorithmic}
\end{algorithm}

All together, the full cost of forming $\t{m}_k$ \eqref{eq_full_model_sketched} is $\bigO(d^3 + md^2 + nd)$ flops for sampling and hashing or $\bigO(d^3 + md^2 + mnd)$ for Gaussian sketching, compared to $\bigO(nd^2)$ flops for constructing the non-sketched model $m_k$ \eqref{eq_full_model}.
In our experiments, we use $m=\bigO(d)$, so the sketching cost becomes $\bigO(d^3 + nd)$ for sampling and hashing or $\bigO(d^3 + nd^2)$ for Gaussian sketching, where again $n\geq d$.
The costs of all other algorithm components are unaffected by the use of sketching, but, as shown in Figure \ref{fig_dfols_times}, these are not the dominant costs in practice, nor are they dominant asymptotically \citep[Table 7.1]{Roberts2019}.

Hence the computational cost is always linear in the data dimension $n$, but, in the big data regime ($n\gg d$), Algorithm \ref{alg_dfols_sketching} with sampling/hashing has linear cost in $d$, rather than quadratic without sketching (i.e.~$\bigO(nd)$ rather than $\bigO(nd^2)$).
Thus we expect Algorithm \ref{alg_dfols_sketching} to have faster runtime than not using sketching, and this to be most significant for large $d$ (provided the ratio $n/d\gg 1$ remains fixed). 

\section{Numerical Results} \label{sec_numerics}
To test Algorithm \ref{alg_dfols_sketching}, we modify DFO-LS to use the three sketching variants described above.
We compare these with the original DFO-LS with by testing on two overdetermined and large-scale CUTEst problems \cite{Gould2015}. 
DFO-LS was shown to have strong performance compared to other solvers in \cite{Cartis2018}.
All problems are run for a maximum of $2(d+1)$ objective evaluations\footnote{Comparable to 2 Jacobian evaluations using finite differences.} and a runtime of 4 hours.
Below, we plot the best objective value achieved after a given runtime for each variant (averaging 10 independent runs for the sketching variants). 
\vspace{-0.5em}


\begin{figure}[tbh]
    \vskip 0.0in
    \begin{center}
        \centerline{\includegraphics[width=\columnwidth]{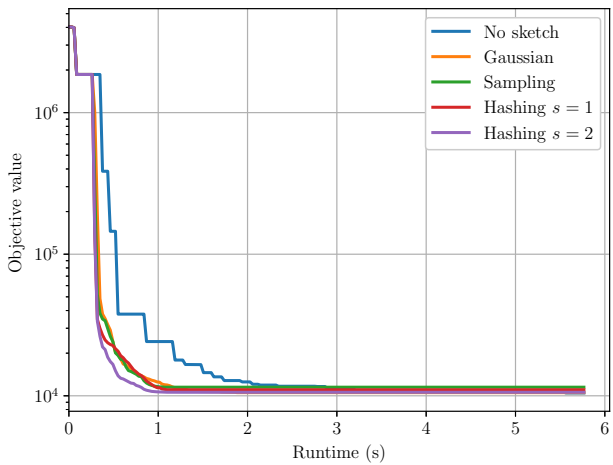}}
        \caption{Objective decrease vs.~runtime for problem DMN15103. Comparing different sketch types, all with $m=d$.}
        \label{fig_dmn15103_compare_method_1n}
    \end{center}
    \vskip -0.2in
\end{figure}

\vspace{-0.3em}
First, in Figure \ref{fig_dmn15103_compare_method_1n} we compare the objective reduction by runtime of the different sketching methods, all with $m=d$, for problem DMN15103 ($d=99$, $n=4643$).
Overall, we find that all sketching variants decrease the objective faster than the original (no sketching) version of DFO-LS and reach the same objective value.
Of the different variants, hashing with $s=2$ achieves the fastest objective decrease. 
Hashing with $s=1$, sampling and Gaussian sketching achieve very similar results.



\begin{figure}[tbh]
    \vskip 0.1in
    \begin{center}
        \centerline{\includegraphics[width=\columnwidth]{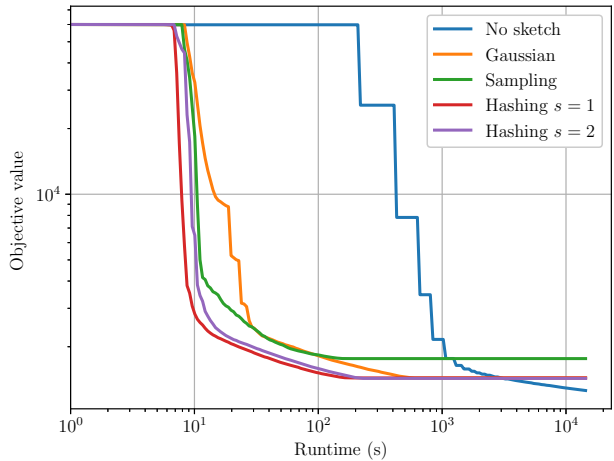}}
        \caption{Objective decrease vs.~runtime for problem MNIST0. Comparing different sketch types, all with $m=d$.}
        \label{fig_mnist0_compare_method_1n}
    \end{center}
    \vskip -0.2in
\end{figure}

We now consider a larger problem, MNISTS0 ($d=494$, $n=60000$), corresponding to training a logistic classifier for zero digits in MNIST \cite{Boyd2018}.
In Figure \ref{fig_mnist0_compare_method_1n} we compare the objective decrease by runtime for different sketch types (using $m=d$).
The `no sketch' variant used the full runtime, and all sketching variants terminated on the evaluation budget.
Here, `no sketching' achieves the lowest objective value, followed by Gaussian sketching and hashing, but all sketching variants have runtime at least one order of magnitude lower.
That is, we gain dramatically on runtime with a small loss of robustness.

\begin{figure}[H]
    \vskip 0.2in
    \begin{center}
        \centerline{\includegraphics[width=\columnwidth]{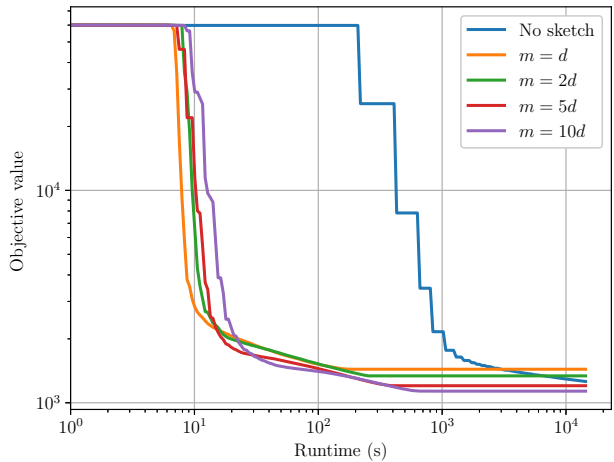}}
        \caption{Objective decrease vs.~runtime for problem MNIST0. Comparing different values of $m$ for hashing with $s=1$.}
        \label{fig_mnist0_compare_size_hashing1}
    \end{center}
    \vskip -0.2in
\end{figure}

This small loss of robustness can be addressed by choosing a larger $m$.
In Figure \ref{fig_mnist0_compare_size_hashing1} we show the results of hashing with $s=1$ for different choices of $m$.
For all $m$, we still gain an order-of-magnitude runtime decrease, but we also reach a comparable (or better) objective value than the original DFO-LS for $m\geq 5d$, in this low-accuracy regime. 

\section{Conclusions and Future Work}
We propose using sketching to reduce the model construction cost---the dominant computational cost---of a DFO algorithm for nonlinear least-squares problems.
With sparse sketches, this reduces the per-iteration computational cost with $n$ residuals and dimension $d$ from $\bigO(nd^2)$ to $\bigO(nd)$ in the big data regime ($n\gg d$).
On large-scale problems, we get the same objective decrease as state-of-the-art software in the low-accuracy regime, but with an order-of-magnitude runtime decrease.
Directions for future work include applying Algorithm \ref{alg_dfols_sketching} to training deep neural networks, developing convergence theory, and using sketched residuals to measure progress \eqref{eq_rho}, similar to \cite{Cartis2018a}.
More detailed experiments would help to select the sketch size $m$ and type (including varying $m$ between iterations, which may aid in the high-accuracy regime).



\bibliography{sketching_refs}
\bibliographystyle{icml2020}

\end{document}